\input epsf

\input amssym.def

\magnification=1200

\baselineskip=14pt
\def\qed{$\vrule height4pt depth0pt width4pt$}
\def\ss{\smallskip}
\def\ms{\medskip}
\def\bs{\bigskip}
\def\D{\Delta}
\def\Z{\bf Z}
\def\ni{\noindent}

\def\<{\langle\langle}
\def\>{\rangle\rangle}

\def\l{\lambda}
\def\g{\gamma}

\line{}\vskip -.5truein
%March 2004%
\centerline{\bf Knot Group Epimorphisms}
\ms
\centerline{DANIEL S. SILVER and WILBUR WHITTEN}
 \bs

\noindent {\narrower\narrower\smallskip\noindent {\bf
Abstract:}     Let $G$ be a finitely generated group, and let $\lambda \in G$. If there exists a knot $k$ such that $\pi k = \pi_1(S^3\setminus k)$ can be mapped onto $G$ sending the longitude to $\lambda$, then there exists infinitely many distinct prime knots with the property. Consequently,  if $\pi k$ is the group of any knot (possibly composite), then there exists an infinite number of prime knots $k_1, k_2, \cdots$ and epimorphisms $\cdots \to \pi k_2 \to \pi k_1 \to \pi k$ each perserving peripheral structures. Properties of a related partial order on knots are discussed. \bs}

\ni {\bf 1. Introduction.} Suppose that 
$\phi: G_1 \to G_2$ is an epimorphism of knot groups
preserving peripheral structure (see \S2). We are motivated by
the following questions.
\bs

\ni {\bf Question 1.1.} If $G_1$ is the group of a prime knot, can
$G_2$ be other than $G_1$ or  $\Z$? \ms

\ni {\bf Question 1.2.} If $G_2$ can be something else, can it be
the group of a composite knot? \bs

\ni Since the group of a composite  knot is an  amalgamated  product of 
the groups of the factor knots, one might expect the answer to 
Question 1.1 to be no.  Surprisingly, the answer to
both questions is yes, as we will see in \S2. 

These considerations suggest a natural partial ordering on knots:
$k_1
\ge k_2$ if the group of $k_1$ maps onto the group of $k_2$
preserving peripheral structure. We study the relation in \S3. 
\bs

\ms
\footnote{} {First author is partially supported by NSF grant
DMS-9704399.}
\footnote{}{2000 {\it Mathematics Subject Classification.}  
Primary 57M25; secondary 37B10.}

\ni {\bf 2. Main result.}  As usual a  knot is the image of
a smooth embedding of a circle in $S^3$. Two knots
are equivalent if they have the same {\it knot type}, that is, there
exists an autohomeomorphism of $S^3$ taking one knot to the other. 

Let
$k$ be a knot in
$S^3$. We denote its group $\pi_1(S^3 \setminus {\rm int}V, *)$
by $\pi k$. Here $V \cong k \times D^2$ is a tubular
neighborhood of $k$, and the basepoint $*$ is chosen on the
boundary $\partial V\cong k \times S^1$. The element  $m$ represented by
a simple closed curve in $\partial V$ that is contractible in $V$ is called a {\it meridian}; the element $l$
represented  by a simple closed curve in $\partial V$ that is
nullhomologous in $S^3
\setminus {\rm int}V$ is  called a {\it longitude}. A well-known
algorithm enables one to express $l$ in terms of Wirtinger
generators corresponding to a diagram for $k$. Details can be found
on page 37 of [BZ85].

The inclusion map $ \partial V \hookrightarrow S^3 \setminus {\rm
int}V$ induces an  injection of fundamental groups. Its image
is the subgroup $\langle m, l\rangle$ generated by  $m$ and $l$. 

Let $k_i, i= 1,2$, be knots with meridian-longitude pairs $m_i,
l_i$. A homomorphism
$\pi k_1 \to \pi k_2$ {\it preserves peripheral structure} if the
image of $\langle m_1, l_1\rangle$ is conjugate to a subgroup of
$\langle m_2,  l_2\rangle$. 
\bs

\ni {\bf Definition 2.1.} (i) $k_1$ {\it covers} 
$k_2$ (or $k_2$ {\it supports} $k_1$) if there is an
epimorphism
$\pi k_1 \to \pi k_2$;

(ii) If $G$ is a finitely generated group normally generated by an
element $\mu$, and if $\l \in G$, then a knot $k$ {\it covers}
$(G, \mu, \l)$ (or briefly $k$ {\it covers} $G$ ) if there
exists an epimorphism $\phi: (\pi k, m, l) \to (G, \mu, \l)$,
where $(m,l)$ is a meridian-longitude pair.\bs

If $k$ covers
$G$ for  a given $\l \in G$, then we say (after Johnson and
Livingston [JL89]) that $k$ {\it realizes} $\l$.
For a given group $G$ as above and $\l\in G$, 
[JL89] provides necessary and sufficient conditions for the
existence of a knot $k$ that covers $(G, \mu, \l)$. We will show
that $k$ can always be chosen to be a prime knot. \bs

\ni {\bf Theorem 2.2.} Let $G$ be a finitely generated group that is
normally generated by a single element $\mu$, and let $\l\in G$.
If there exists a knot $k$ that realizes $\l$, then there exists
an infinite number of distinct prime knots
that  realize  $\l$. \bs

\ni {\bf Proof.} As the theorem easily follows when $k$ is
trivial, we assume that $k$ is knotted. By Proposition 2.5 of
[EKT03], we can regard
$k$ as the numerator closure $T^N$ of a prime tangle $T$.  (See  [EKT03], where   the authors provide general terminology and prove an even stronger condition). Set
$C = (1/2)^N$, using Conway notation (Figure 1). \bs

%\centerline{\includegraphics[width=.5in]{KG1.pdf}} \bs
\epsfxsize=.5truein
\centerline{\epsfbox{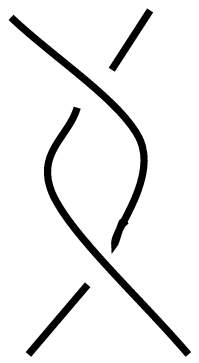}} \bs

\centerline{{\bf Figure 1:} The tangle $1/2$}\bs

We use a construction of [EKT03] to form a 2-component link 
$L= T^N \cup (1/2)^N  = k \cup C$ (Figure 2). The components
$k$ and $C$ are contained in disjoint solid tori $V_1$ and 
$V_2$, respectively, the cores of which form  a Hopf link. Note
that $C$ is an untwisted double of the core of $V_2$. 
By Propositions 2.3 and 2.4 of [EKT03],  $L =k\cup C$ is a prime
link.
\bs

%\centerline{\includegraphics[width=2in]{KG2.pdf}} \bs
\epsfxsize=2truein
\centerline{\epsfbox{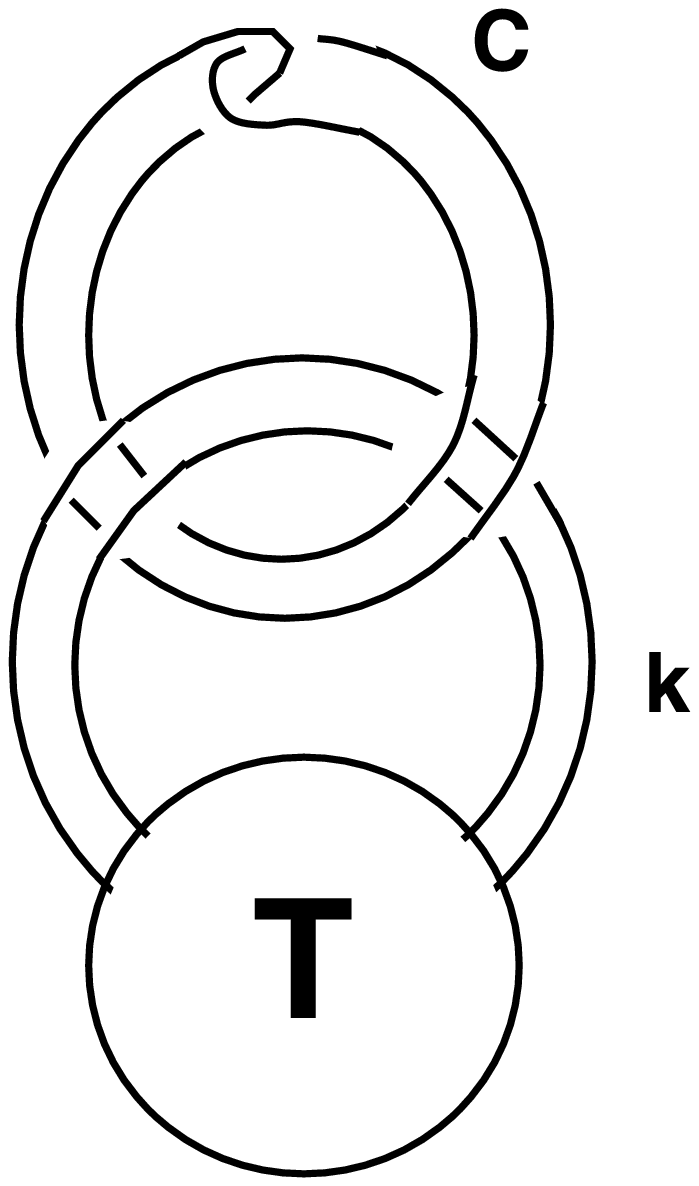}} \bs

\centerline{\bf Figure 2: The link $L$}\bs

Let $K_q$ denote the knot obtained from $k$ by $1/q$-surgery on
the trivial knot $C$. We take $q>0$ for simplicity. It is helpful to view $K_q$ as  a satellite
of a $(-2q)$-twist knot $\g_{2q}$. (The notation $\g_{2q}$ here means
that the usual regular projection of $\g_{2q}$ has $-2q$ negative
half-twists.) From this perspective, it follows immediately that each $K_q$ is nontrivial. 

If $k$ is
composite, then it follows from Theorem 4.1 of [G-S97]
that $K_q$ is prime for $q\ge 3$, since $k\cup C$ is prime. On
the other hand, if $k$ is prime and $K_{q'}$ is composite
for some $q'$, then $K_q$ will be prime when $q\ge q'+3$.
Thus there exists $d\ge 3$ such that $K_q$ is prime whenever $q \ge d$;
cf. Remarks 3.5 (ii).

For a fixed meridian-longitude pair $(m, l)$ of $k$, we have an
epimorphism $\phi: \pi k \to G$ for which $\phi(m)=\mu$ and
$\phi(l)=\l$. The elements $m$ and $l$ are represented by oriented simple
closed curves $m_1$ and $l_1$ on the boundary of a tubular neighborhood
of $k$ in $V_1$. Assume that $m_1 \cap l_1$ is a single point $*$,
which we take as the basepoint for $\pi(k \cup C), \pi K_q$ and 
$\pi k$, as previously mentioned. The curve $m_1$ certainly represents a meridian
of $K_q$, and since the linking number of $k$ and $C$ is zero, the curve
$l_1$ represents a longitude (see pages 266--267 of [R76]). Furthermore,
if $(w,c)$ is a meridian-longitude pair for $C$ (based at $*$), we 
have $\pi k = \pi(k \cup C)/\< w \>$ and $\pi K_q =
\pi(k \cup C)/\< wc^q\>$. (Here and throughout $\<\ \>$ denotes normal
closure.) Thus
$$\pi K_q/\<w\> = \pi(k \cup C)/\< w,
wc^q\> = \pi k/\< c^q\>.$$
But $c=1$ in $\pi k$, since $C$ represents the identity element of 
$\pi k$, and so $\pi K_q/\< w \> = \pi k$. The induced
homomorphism $\eta: \pi K_q \to \pi k$ is an epimorphism taking
meridian to meridian and longitude to longitude, preserving orientations; cf. Remarks 3.5 (ii).
Therefore, $\phi_q = \phi \circ \eta: \pi K_q \to G$ is an epimorphism
that maps the classes of $m_1$ and $l_1$ to $\mu$ and $\l$,
respectively.

Recall that $K_q$ ($q \ge d$ from now on) is a satellite of a $(-2q)$-twist knot  $\g_{2q}$. Since $\g_{2q}$ is  the double of the unknot with 
twisting number $-2q$, it follows from the form of the Alexander polynomial for doubled knots (see p. 136 of [BZ85], for example)  that the exteriors
${\rm Ext}(\g_{2q})$ and ${\rm Ext}(\g_{2\bar q})$ are not
homeomorphic if $q\ne \bar q$. 

Finally, consider the canonical splitting of the exterior $E_q$ of
$K_q$ into a union $\Sigma_q$ of Seifert pieces (the characteristic
submanifold of $E_q$) and a union $\Lambda_q$ of atoroidal pieces
($\Lambda_q = {\rm cl}(E_q \setminus \Sigma_q)$) [JS79], [J76]. Each
of 
$\Sigma_q$ and $\Lambda_q$ has a finite number of components, and
$${\rm cl}[(\Sigma_q \cup \Lambda_q) \setminus {\rm Ext}(\g_{2q})] 
\cong {\rm cl}[(\Sigma_{r}\cup \Lambda_{r}) \setminus {\rm
Ext}(\g_{2r})],$$
for $q,r \ge d$. Now choose  $p \in \{d,d+1,\ldots \}$ so large
that a copy of ${\rm Ext}(\g_{2q})$ is not a component of
$\Lambda_q \setminus {\rm Ext}(\g_{2q})$, for any $q\ge p$. Thus
if $q\ne r$ and $q, r \ge p$, then $E_q$ is not homeomorphic
to $E_{r}$, since ${\rm Ext}(\g_{2q})$ is not homeomorphic to 
${\rm Ext}(\g_{2r})$. Therefore, there exist infinitely many
distinct prime knots realizing $\l$. \qed\bs

\ni {\bf Remarks 2.3.} (i) Notice the seemingly large number of choices 
we have in the  above construction of prime knots that realize
$\l$. For example, we might double and redouble $C$ itself. \ss

(ii) To see that the answer to each of Questions 1.1 and 1.2 is yes, see
Example 2.6. Less specifically, let $G$ be the group of a composite knot
$k$ with meridian-longitude pair
$(\mu, \l)$, and let $k_1$ and $k_2$ be ambient isotopic copies
of $k$. Assume that $k$ is oriented and that each of $k_1$ and
$k_2$ inherits this orientation. Let $(m_i, l_i)$ be a
meridian-longitude pair for $k_i$ ($i = 1,2$), and let $\phi_i: \pi
k_i \to \pi k$ be an isomorphism  such that $\phi(m_i) = \mu$ and
$\phi_i(l_i)= \l$. Then $\phi_1$ and $\phi_2$ induce an epimorphism
$\phi: \pi(k_1\sharp k_2) \to G$ such that $\phi(m) = \mu$ and
$\phi(l) = \l^2$, for some choice of meridian-longitude pair
$(m,l)$ for $k_1\sharp k_2$. Thus $k_1\sharp k_2$  covers $(G, \mu,
\l^2)$, and so there exist prime knots $K_q$ covering $(G, \mu,
\l^2)$ by Theorem 2.2.

 \bs

\ni {\bf Corollary 2.4.} Let $G$ be a knot group normally generated
by $\mu \in G$, and let $\l \in G$. Then there exists an infinite
number of prime knots realizing $\l$ if and only if $\l \in G'' \cap
Z(\mu)$, where $G''$ is the second commutator subgroup of $G$ and 
$Z(\mu)$ is the centralizer of $\mu$ in $G$. \bs

\ni {\bf Proof.} According to the main result, Proposition 1, of 
[JL89], $\l$ is realizable (by some knot) if and only if $\l \in 
G'' \cap Z(\mu)$, since $G$ is a knot group. But Theorem 2.2 implies that
$\l$ is realizable if and only if it is realizable by an infinite number of
prime knots.
\qed\bs

\ni {\bf Remarks 2.5.} (i) If we take $\l$ to be the longitude of a knot
group $G_2$ with meridian $\mu$, then Corollary 2.4 ensures the
existence of a knot group $G_1$ not isomorphic to $G_2$  and an
epimorphism
$\phi: G_1
\to G_2$ that preserves peripheral structure. This provides another
answer to Questions 1.1 and 1.2. 

(ii) One can avoid Corollary 2.4 by noting that the group of any knot $K$
covers itself by the identity automorphism $(\pi K, m, l) \to (\pi
K, m, l)$, and then taking $k = K$ in the construction given in
the proof of Theorem 2.2. \bs

\ni {\bf Example 2.6.} Consider the  2-component link $L = k \cup C$ in Figure 3. Regard the knotted component $k$ as the connected sum $(k_W \sharp k_N) \sharp (k_E \sharp k_S)$, of two Granny knots. (Here $W, N, E, S$ abbreviate west, north, east and south.) As in Remarks 2.3 (ii), there exist isomorphisms $\pi (k_W \sharp k_N) \to \pi (k_W \sharp k_E)$ and $\pi (k_E \sharp k_S) \to \pi (k_W \sharp k_E)$, taking the longitude of each of $\pi (k_W \sharp k_N)$ and $\pi (k_E \sharp k_S)$ to that of $\pi (k_W \sharp k_E)$. Hence there is an epimorphism  $\phi: \pi k \to \pi (k_W \sharp k_E)$ that takes a longitude of $\pi k$ to the square of that of $\pi  (k_W \sharp k_E)$.  Moreover,  $\phi$ maps the class of  $C$ trivially. 

Let $K_q$ be the knot resulting from $k$ after $1/q$ surgery on $C$. As in the proof of Theorem 2.2, the group $\pi K_q$ admits an epimorphism onto the Granny knot group, sending meridian to meridan, longitude to the square of a longitude

We obtain the conclusion of Theorem 2.2 by another method, one that enables us to obtain hyperbolic knots $K_q$. Since the link $L$ is prime and alternating, a theorem of W. Menasco  [Me84] implies that $L$ is hyperbolic.  Results of W. Thurston [Th77/83]  and W. Neumann and D. Zagier [NZ85] imply that for sufficiently large $q$ the knots $K_q$ are hyperbolic, with strictly increasing volumes (that approach 36.4732$\ldots$); in particular, the knots are distinct. \bs

%\centerline{\includegraphics[width=2in]{KG3.pdf}} \bs
\epsfxsize=2truein
\centerline{\epsfbox{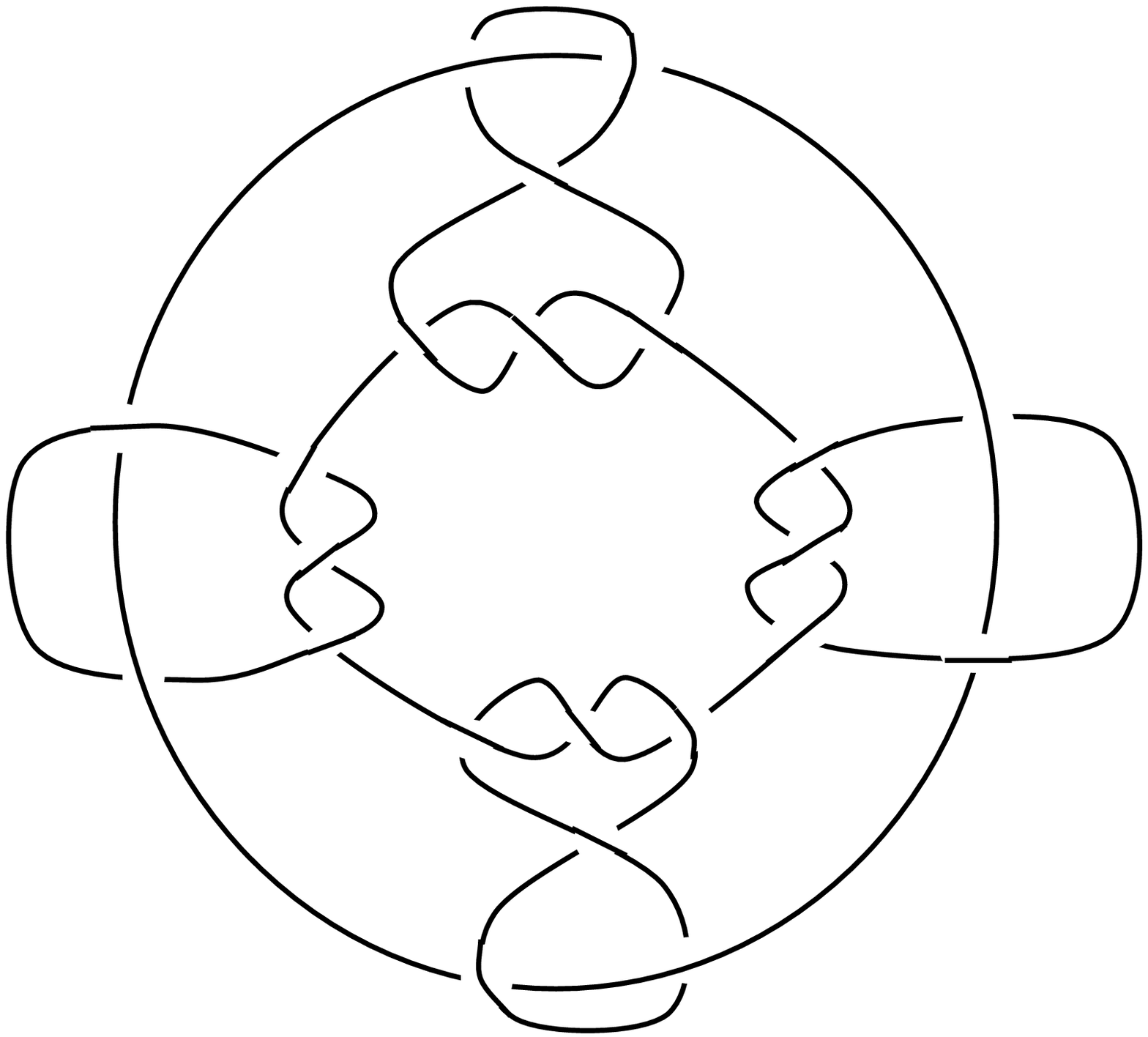}} \bs

\centerline{\bf Figure 3:  $k \cup C$}\bs

\ni {\bf Remark 2.7.} The choice of a curve $C$ as in Example 2.6 and in Theorem 2.2 satisfies four basic requirements: \ss

\item{} 1. $C$ is unknotted in $S^3$; 
\item{} 2. The linking number of $C$ and $k_1 \sharp k_2$ is zero; 
\item{} 3. The link $C \cup (k_1 \sharp k_2)$ is prime;  and 
\item{} 4. The epimorphism $\phi: \pi (k_1 \sharp k_2) \to \pi k$ maps the class of $C$  trivially. \bs

\ni {\bf 3. Partial order on knots.} Motivated by the 
results of \S2, we introduce a relation on the set of all knots. \bs

\ni {\bf Definition 3.1.} $k_1 \succeq k_2$ if there is an epimorphism
$\pi k_1 \to \pi k_2$ preserving peripheral structure.  \bs

\ni {\bf Proposition 3.2.} The relation $\succeq$ is a partial order.\bs

\ni {\bf Proof.} Clearly $\succeq$ is reflexive and
transitive. It remains to show that $\succeq$ is antisymmetric.

If
$k_1 \succeq k_2$ and $k_2 \succeq k_1$, then there exist epimorphisms
$\phi: \pi k_1 \to \pi k_2$ and $\eta: \pi k_2 \to \pi k_1$,
preserving peripheral structure. The compositions $\phi\circ \eta$
and $\eta\circ \phi$ are epimorphisms from a knot group to itself.
Since any knot group is residually finite
[H87] and finitely generated,  it has the Hopfian property:
any epimorphism from the group to itself is an isomorphism [M40]. Hence
both
$\phi\circ
\eta$ and $\eta\circ \phi$ are isomorphisms. In particular, it follows
that 
$\phi$ is an isomorphism preserving peripheral structure. By [W68]
$k_1$ and $k_2$ have homeomorphic complements. Finally [GL89]
implies that $k_1$  and $k_2$ are of the same knot type (not necessarily of the same ambient isotopy type, however, since chiral knots exist).\qed\bs

\ni {\bf Remarks 3.3.} (i) The condition that the epimorphism preserve
peripheral structure is needed for the conclusion of Proposition
3.2. To see that this is so, let $k_1$ be the granny knot and 
$k_2$ the square knot. Since their groups are isomorphic, we 
certainly have epimorphisms from one  to the other. However, no 
such epimorphism preserves peripheral structure, and indeed 
$k_1$ is not equal to $k_2$. \ss

(ii) The
relation
$\succeq$ is compatible with some well-known invariants. For example, if
$\D_k^{(i)}(t)\ (i \ge 1)$ denotes the $i$th Alexander polynomial of $k$, then 
$k_1 \succeq k_2$ implies the existence of an epimorphism $\pi k_1 \to \pi
k_2$, which in turn implies that, for each $i$, 
$\D^{(i)}_{k_1}(t)$ contains
$\D^{(i)}_{k_2}(t)$ as a factor. Following the usual practice,  we will abbreviate $\D^{(1)}_k(t)$ by 
$\D_k(t)$ and refer to it as {\sl the} Alexander polynomial of $k$.  Necessary algebraic background information can be found in [MKS76]). \ss

(iii) If $k= k_1 \sharp k_2$,  then $k\succeq
 k_1$ and
$k\succeq k_2$. Hence $\succeq$ refines the crude partial order induced
by knot factorization. \bs
 
There are several natural methods to produce knots $k_1, k_2$ with $k_1 \succeq k_2$. For example, suppose that a diagram for a knot $k_1$ displays a rotational symmetry $f$. Let $k_2$ denote the quotient knot in the $3$-sphere $S^3/f$. By identifying all pairs of Wirtinger generators corresponding to arcs in $f$-orbits we obtain a projection from $\pi k_1$ to $\pi k_2$ preserving peripheral structures. Hence $k \succeq \bar k$. 

Another method is contained in the following. \bs

\ni {\bf Proposition 3.4.} Assume that $k$ is a satellite knot with 
pattern knot $\tilde k$. Then $k\succeq \tilde k$.\bs

\ni {\bf Proof.} The satellite knot $k$ is the image of a diffeomorphism
$g: \tilde V \to \hat V\subset S^3$, where $\tilde V$ is a standard solid
torus containing
$\tilde k$,  and $\hat V$ is a tubular neighborhood of a knot $\hat
k$, the {\it companion knot}. As usual, we require that $\tilde k$ not
be contained in any $3$-ball of $\tilde V$, and also that $g$ send the
longitude $l$ of
$\tilde V$ to the longitude
$\hat l$ of the solid torus
$\hat V$. Denote the meridian of $\hat V$ by $\hat m$.

The group $\pi k$ of
the satellite is isomorphic to the free product of $\pi \hat k$ and
$\pi_1(\tilde V\setminus \tilde k)$ with amalgamation: subgroup $\langle
\hat m, \hat l\rangle$ of $\pi
\hat k$ is identified with the subgroup $\langle  m, l\rangle$ of
$\pi_1(\tilde V\setminus
\tilde k)$, matching $\hat m$ and $\hat l$ with $m$ and $l$ (see [BZ85],
for example). By a well-known property of free products with amalgamation,
both $\pi \hat k$ and $ \pi_1(\tilde V\setminus \tilde k)$ are subgroups
of $\pi k$ embedded in the obvious way. (see [LS77]).

Consider the natural projection $\phi: \pi k \to \pi k/N$ where $N$ is
the normal closure of the commutator subgroup of $\pi\hat k$. Regard $\pi
k/N$ as the quotient group of $\pi k$ obtained by allowing the elements of
$\pi \hat k$ to commute. In the quotient,
$\pi \hat k$ collapses to the infinite cyclic group
generated by $\hat m$.  The effect on
$\pi_1(\tilde V\setminus
\tilde k)$ is to kill $l$, resulting in the group $\pi \tilde k$. Hence
$\phi$ is an epimorphism from $\pi k$ to $\pi \tilde k$. It is obvious
that $\phi$ preserves peripheral structure. \qed \bs

\ni {\bf Remarks 3.5.} (i) Proposition 3.4 states that a satellite
knot
$k$ covers its pattern knot $\tilde k$. It is not generally true that
$k$ covers its companion $\hat k$. To see this, consider the untwisted
double $k$ of the trefoil knot. The Alexander polynomial
$\D_k(t)$ is  trivial (see [BZ85]). However, the companion knot is the
trefoil, which has nontrivial Alexander polynomial. In view of Remark 3.3
(ii), $k$ does not cover $\hat k$.
We mention that for a given satellite knot, the pattern knot might well
cover the companion or vice versa. \ss

(ii) We can avoid the use of Theorem 4.1 of [G-S97] and apply
Proposition 3.4 to give an alternative proof of Theorem 2.2 as follows. Think of
$K_q$ as the satellite knot with pattern $k$ in the interior of the
(standardly embedded) solid torus $V_1$ and with companion the twist-knot
$\gamma_{2q}$; here $V_1$ is mapped by a longitude-preserving homeomorphism
onto a tubular  neighborhood of $\gamma_{2q}$. Since ${\rm Int}(V_1)$
contains no $2$-sphere that decomposes $K_q$ as a nontrivial connected sum,
it follows easily that $K_q$ is prime for $q \ge 2$, say. Proposition 3.4
immediately yields the epimorphism $\eta: \pi K_q \to \pi k$. We omit
details. \bs

%Proposition 3.3 below provides another example. 

%In [G82] M. Gromov introduced a simplicial volume invariant
%$||M||_0$ for closed manifolds $M$. The invariant was extended for 
%compact $3$-manifolds with toroidal boundary by W. Thurston [T78],
%who proved that if $M$ is hyperbolic then $||M||_0$ is equal to 
%the  hyperbolic volume of $S^3 \setminus k$ divided by the volume
%of a regular ideal simplex. He also showed that if there exists a 
%degree $d$ map from $M_1$ to $M_2$, then 

%$$||M_2||_0 \ge d||M_1||_0 \eqno(3.1)$$

%There is a natural  invariant $||\ell||_0$
%for a link $\ell \subset S^3$, defined as the Gromov volume of its
%exterior. It has been studied by T. Soma in [S81].  
%\bs

%\ni {\bf Proposition 3.4.} If $k_1\ge k_2$ then $||k_1||_0\ge
%||k_2||_0$.
%\bs

%\ni {\bf Proof.} An  epimorphism $\pi k_1 \to \pi k_2$ preserving
%peripheral structure maps the meridian-longitude pair $m_1, l_1$ to
%a elements that are pairwise-conjugate to  $m_2^\epsilon l_2^p,
%l_2^q$, for $\epsilon = \pm 1$ and some integers $p, q$.  

% Since every knot exterior is a
%$K(\pi, 1)$, the condition that 
%$k_1\ge k_2$ is equivalent to the requirement that  there exists
%a degree $d$ map from the exterior of $k_1$ to that
%of $k_2$. The Proposition follows now from (3.1) above. \qed\bs

%In [R92] Y. Rong studied the partial order on the set of compact
%$3$-manifolds with toroidal boundary: $M_1\ge M_2$ if there exists
%a degree one map from $M_1$ to $M_2$. The following is an 
%immediate consequence of his main result, proved with the help of
%Gromov's norm. 
%\bs

Recall that the {\it genus} $g(k)$ of a knot $k$ is the smallest
genus of any Seifert surface of $k$. \bs

\ni {\bf Conjecture 3.6.} If $k_1\succeq k_2$ then $g(k_1) \succeq g(k_2)$.
\bs

As evidence for Conjecture 3.6 we offer the following. \bs

\ni {\bf Proposition 3.7.} Assume that $k_1 \succeq k_2$. If either
(i) $k_1$ is fibered or (ii) $k_2$ is alternating then $g(k_1) \succeq
g(k_2)$. \bs

\ni {\bf Proof.} (i) A knot $k$ is fibered if and only if the
commutator subgroup $[\pi k, \pi k]$ is finitely generated, in
which case it is a free group of rank equal to $2g(k)$ (see
[R76], for example). Assume that
$k_1 \succeq k_2$. The epimorphism  $\pi k_1 \to \pi k_2$ restricts to an
epimorphism of commutator subgroups. Consequently, $[\pi k_2, \pi
k_2]$ is finitely generated, and hence free of rank less than
or equal that  of $[\pi k_1, \pi k_1]$. Since the rank of $[\pi
k_2, \pi k_2]$ is equal to $2g(k_2)$, the proof is
complete. 

(ii) By a theorem of  Seifert [S34], for any knot $k$ we have
$2 g(k) \ge {\rm deg} \D_k(t)$, while a theorem of Murasugi
[M60] states that equality holds when 
$k$ is alternating. Assume that
$k_1 \succeq k_2$. Then $\D_{k_2}(t)$ divides $\D_{k_1}(t)$.
Hence $2g(k_1) \ge  {\deg} \D_{k_1}(t) \ge {\rm
deg}\D_{k_2}(t) = 2g(k_2)$. \qed \bs

There is further evidence for Conjecture 3.6. If $k_1 \succeq k_2$, then the epimorphism
$\pi(k_1) \to \pi(k_2)$ is induced by a boundary-preserving map $f: {\rm Ext}(k_1) \to {\rm Ext}(k_2)$. If $S$ is a minimal genus Seifert surface spanning $k_1$, and if $f_*([S]) \in H_2({\rm Ext}(k_2), \partial {\rm Ext}(k_2)) (\cong {\Bbb Z})$ is not zero, then it follows from Corollary 6.22 of [G83] that
$g(k_1) \ge g(k_2)$. We are grateful to Ian Agol for pointing out the connection between Conjecture 3.6 and Gabai's result. \bs

\ni {\bf Proposition 3.8.} If $k_i$ is a sequence of hyperbolic
knots such that  $k_0\succeq k_1 \succeq\cdots \succeq k_i \succeq \cdots$, then
$k_i = k_{i+1}$ for sufficiently large $i$. \bs

\ni {\bf Proof.} Assume that $k_0 \succeq k_1 \succeq k_2 \succeq \cdots$,
for hyperbolic knots $k_i$. Consequently, there is a sequence 
of epimorphisms 
$$\pi k_0\ {\buildrel \phi_0 \over \longrightarrow}\ \pi k_1\
{\buildrel \phi_1 \over \longrightarrow}\ \pi k_2\ {\buildrel
\phi_2 \over \longrightarrow}\cdots,$$
each preserving peripheral structure.
By Theorem 1 of [S02], $\phi_i$ is an  isomorphism for  sufficiently
large $i$. (The result of [S02]  requires only knot group epimorphisms without any constraint on peripheral  subgroups.) As in the proof of Proposition 3.2, [W68] and [GL89]
together imply that $k_i = k_{i+1}$ for sufficiently large
$i$.\qed \bs  

\ni {\bf Conjecture 3.9.} (Cf. J. Simon: Problem 1.12 (D) in [K95]) ] 
 Any knot covers only finitely many
knots. In other words, if $k_1$ is a knot, then the collection of all
knots $k_2$ such that $k_1\succeq k_2$ is finite. \bs

Conjecture 3.9 is true if $k_1$ is fibered. To see this, note that 
$k_2$  must also be fibered. In [S95] an entropy invariant $h_k$ was
defined for any fibered knot $k$.  If $k_1 \succeq k_2$, then
$h_{k_1} \ge h_{k_2}$ and also $g(k_1)\ge g(k_2)$. 
By Theorem 3.4 of [S95], there exist only finitely many fibered
knots with genus and entropy no greater than given bounds. \bs

\ni {\bf Definition 3.10.} A knot $k$ is {\it minimal} (with respect
to the partial order) if $k \succeq k'$ implies that $k= k'$
or else $k'$ is trivial.  \bs

\ni {\bf Proposition 3.11.} If $k$ is a fibered knot with irreducible
Alexander polynomial, then $k$ is minimal. \bs

\ni {\bf Proof.} Assume that $k\succeq k'$. Since $\D_{k'}(t)$ divides
$\D_k(t)$ and $\D_k(t)$ is irreducible, either $\D_{k'}(t) =
\D_k(t)$ or else $\D_{k'}(t) = 1$. In the first case, the epimorphism
$\phi: \pi k \to \pi k'$ restricts to an isomorphism  of commutator
subgroups, as these groups are both free of rank equal to ${\rm deg}\D_k(t)$; it follows that $\phi$ is itself an
isomorphism, and as in the proof of Proposition 3.2, $k= k'$. 
In the second case, $\D_{k'}(t) =1$, and this together with the fact
that $k'$ is fibered imply that  $k'$ is trivial. \qed \bs

The figure eight knot $4_1$ is minimal by
Proposition 3.11. However, Example 3.12 below shows that $4_1$
does not remain minimal in the  larger category of virtual knots. 
It shows also that Conjecture 3.9 is  not  true in the category.

A knot
is often studied as an equivalence class of planar knot
diagrams, two diagrams being equivalent if one can be obtained from
the other by a sequence of Reidemeister moves. In 1997 Kauffman
introduced  virtual knot diagrams, allowing a new type of crossing,
called a {\it virtual crossing} and 
indicated by a small circle surrounding the site. After
suitably extending the usual
Reidemeister moves to allow certain deformations involving virtual
crossings, Kauffman defined a {\it virtual knot} to be an 
equivalence class of virtual diagrams. The reader is referred to
[K97], [K99], [K00] for details. It is a remarkable feature of  the
theory that two classical knots are equivalent under  generalized
Reidemeister moves if and only if they are equivalent under the
classical ones [GPV00]. In this sense, virtual knot theory is an
extension of the  classical theory. 

Many classical invariants of knot theory extend naturally in 
the larger virtual category. In particular, one can associate
a knot group to an equivalence class of diagrams. Virtual knot
groups were classified in [SW00] both in terms of combinatorial
presentations  and topologically (see also [Ki00]). The peripheral
structure of a virtual knot is defined just as for classical
knots (see [Ki00] for details). \bs

\ni {\bf Example 3.12.} A diagram for the figure eight knot appears
in Figure 4 with labeled Wirtinger generators. The group of
$k$   has a presentation
$\langle x,y, z, w \mid zx=yz, yw=zy, wx=zw\rangle$. (The
fourth Wirtinger relation $xw=yx$ is a consequence of the other
three, and so we have omitted it.) Using the second and third
relations to express $w$ and $z$ in terms of $x$ and $y$, we 
see that  
$$\pi k \cong \langle x, y \mid x^{-1} y^{-1} x y^{-1} x^{-1} y x
y^{-1} xy\rangle.$$

Consider the diagram for the virtual knot $k_q$ ($q \ge 2$) in
Figure 4 with  labeled Wirtinger generators. The corresponding group
generated by
$x,y,z,w, w_1, \ldots, w_{q-1}$ with relations $zx=yz, yw=zy,
wx=zw, xw_1 = wx, xw_2= w_1x, \ldots, xw_{q-2}=w_{q-3}x,
wy=w_{q-2}x$.  \bs

%\centerline{\includegraphics[width=2.2in]{KG4.pdf}} \bs
\epsfxsize=2.2truein
\centerline{\epsfbox{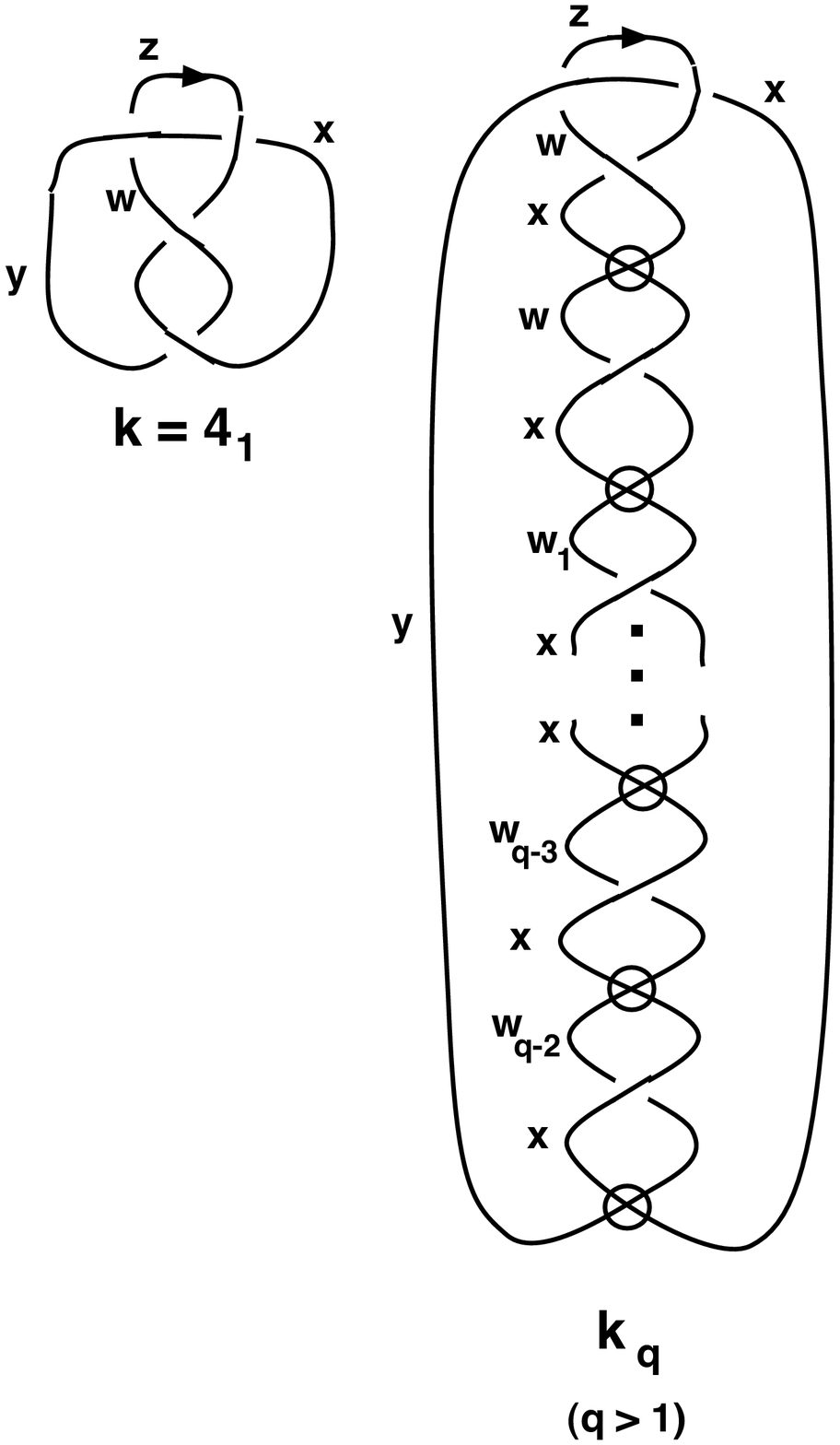}} \bs

\centerline{\bf Figure 4:  Figure eight knot and $K_q$}\bs

The fourth relation implies that
$w_1=x^{-1}wx$. The fifth implies that $w_2=x^{-2}wx^2$,
and so forth. The next to last relation implies that
$w_{q-2} = x^{-(q-2)} w w^{q-2}$. We use these to eliminate
$w_1, \ldots, w_{q-2}$. The last relation then becomes 
$x^{-q} y x^q = y$. Consequently,

$$\pi k_q \cong \langle x,y,z, w\mid zx=yz, yw=zy,
wx=zw, y= x^{-(q-1)}y x^{q-1} = w \rangle.$$
Recall that the three relations
$zx=yz, yw=zy, wx=zw$ together imply that $xw=yx$, which can be
rewritten as $w= x^{-1}yx$. Substitution in the relation 
$y= x^{-(q-1)}y x^{q-1} = w$ yields

$$\pi k_q \cong \langle x,y,z, w\mid zx=yz, yw=zy,
wx=zw, x^{-q} y x^q = y \rangle$$
$$\cong \pi k/\<x^{-q} y x^q = y\>$$
$$\cong \langle x, y \mid x^{-1} y^{-1} x y^{-1} x^{-1} y x
y^{-1} xy, x^{-q} y x^q = y\rangle$$.

The Reidemeister-Schreier method (see [LS77], for example) can be
used to find a presentation for the commutator subgroup of $\pi
k_q$:
$$\langle a_j \mid a_{j+2} = a_{j+1}^2 a_j^{-1} a_{j+1},\ 
a_{j+q}=a_j\quad (j\in \Z)\rangle.$$
(This group is in fact the fundamental group of the $q$-fold cyclic
cover of $S^3$ branched over $k$. However, we do not require this
fact.) Its abelianization is finite, and it has order equal to the
absolute value of the cyclic resultant ${\rm Res}(t^2-3t+1, t^q-1)$.
These values grow exponentially with $q$ (see [SW02]). Consequently,
the groups $\pi k_{q}$ are pairwise nonisomorphic for sufficiently
large $q$. (In fact, they are all nonisomorphic, but this is again a
fact that we do not require.) 

The canonical projection $\phi_q: \pi k \to \pi k_q$ is a surjection that 
maps the elements $x$ and $z^{-1}xy^{-1}w$, forming a peripheral pair 
for $k$, to their cosets in $\pi k_q$. It is easy to see that 
$(x, z^{-1}x^{-(q-1)}y^{-1}wx^{n-2})$ is  a peripheral pair for $k_q$.
Using the relations above, one checks that
$z^{-1}x^{-(q-1)}y^{-1}wx^{n-2}$ is equal to $z^{-1}xy^{-1}w$. Hence the
projection $\phi_q$ preserves peripheral structure, and  we have produced
infinitely many virtual knots $k_q$ such that $k \succeq k_q$.
\bs

\ni {\bf 4. Degree one maps and a related partial order.} 
The partial order in Definition 3.1 is related to another, a 
partial order on compact $3$-manifolds, one that has been studied by Y. Rong [R92], S. Wang [W02], and others:  If $M$ and $N$ are compact, oriented $3$-manifolds,
then $M\succeq_1 N$ if there exists a degree one proper map from $M$ to $N$. (A manifold map
is {\it proper} if it maps boundary to boundary.) 
When
applied to exteriors of knots  $k_1, k_2$, the relation becomes: $k_1 \succeq_1 k_2$ if there is an epimorphism
$\pi k_1 \to \pi k_2$  mapping meridian $m_1$ to $ m_2l_2^p$ and  longitude $l_1$
to $m_2^ql_2^{\pm  1} $, for some integers $p,q$. Since the normal subgroup generated by $ m_2l_2^p$ must be all of $\pi k_2$, Corollary 2 of [CGLS87] implies that $p \in \{0, 1, -1\}$. (The recent proof that every nontrivial knot satisfies Property P [KM04] implies that $p$ must in fact be $0$. However, we do not require that fact.) Also, since $m_2^q l_2^{\pm 1}$ must be in $(\pi k_2)'' \cap Z(m_2)$ [JL89], we have $q=0$.

Clearly, $k_1 \succeq_1 k_2$ implies $k_1 \succeq k_2$.\bs

\ni {\bf Theorem 4.1.} In general, the relation $k_1 \succeq k_2$ does not imply $k_1 \succeq_1 k_2$. \bs

\ni {\bf Proof.} Let $k_1, k_2$ be the torus knots $9_1, 3_1$, respectively.  Since $k_1$ has
a rotational symmetry with $k_2$ as quotient, it follows from the remark before Proposition 3.4
that $k_1 \succeq k_2$. 

Note that the epimorphism $\pi k_1 \to \pi k_2$ induced by the symmetry maps meridian $m_1$ to meridian $m_2$ and longitude $l_1$ to the third power $l_2^3$. It suffices to show there exists no epimorphism mapping $(m_1, l_1)$ to $(m_2l_2^p, l_2^{\pm 1})$, for $p \in \{0, 1, -1\}$. Since torus knots satisfy property P  [S33] (see also [H64]  or  [BZ85], p.274), the integer $p$ must be zero. The desired conclusion follows from the next proposition. (Details are given in Example 4.3.) \qed \bs

The $A$-polynomial was introduced in [CCGLS94]. We briefly review its definition, following
[CL96]. For any knot $k$ with meridian-longitude pair $(m,l)$, consider the affine algebraic set $R = {\rm Hom}(\pi k,  SL_2({\bf C}))$.
Let $R_\bigtriangleup$ be the algebraic subset consisting of all $\rho \in R$ such that $\rho(l)$ and $\rho(m)$ are upper triangular. Let 
$\xi: R_\bigtriangleup \to {\bf C}^2$ be the eigenvalue map $\rho \mapsto (M,L)$, where $M$ and $L$ are the top-left entries (eigenvalues) of  $\rho(m)$ and $\rho(l)$, respectively. The closure of any component of $\xi(R_\bigtriangleup)$ has complex dimension $0$ or $1$. Each $1$-dimensional component is the 
zero set of a polynomial that is unique up to multiplication by a constant; the product of all such polynomials, divided by $L-1$ (the polynomial corresponding to abelian representations) is $A_k$. It is possible to normalize, choosing a 
suitable multiplicative constant, so that the coefficients of $A_k$ are integers with no common divisor, and we do so. \bs

\ni {\bf Proposition 4.2.} Let $k_1, k_2$ be knots. If $\phi: \pi k_1  \to \pi k_2$ is a homomorphism mapping $(m_1, l_1)$ to $(m_2^a l_2^b, m_2^c l_2^d)$, for integers $a,b,c,d$, then each irreducible factor of $A_{k_2}(M,L)$ divides $(M^cL^d -1)\cdot A_{k_1}(M^aL^b, M^cL^d)$. \bs

\ni {\bf Proof.} Denote the eigenvalue maps for $k_1$ and $k_2$ by $\xi_1$ and $\xi_2$, respectively. If $(M,L)$ is in the image of $\xi_2$, then $(M^aL^b, M^cL^d)$
is in the image of $\xi_1$. Hence $A_{k_2}(M,L)  =0$ implies that $(M^cL^d -1)\cdot A_{k_1}(M^aL^b, M^cL^d)=0$.  The result follows from Hilbert's Nullstellensatz (see [H77], for example). \qed \bs

\ni {\bf Example 4.3.} The A-polynomial of  $k_1= 9_1$ is $1+M^{18}L$ while 
that of the trefoil $k_2= 3_1$ is $1+M^6L$ (see [CCGLS94] or [N02]). The 
polynomial $A_{3_1}(M,L)$ divides $A_{9_1}(M, L^3)$, reflecting the fact that 
there is a homomorphism $\pi k_1 \to \pi k_2$ mapping $(m_1, l_1)$ to $(m_2, l_2^3)$. However, the polynomial $A_{3_1}(M,L)$ is irreducible and divides neither
$A_{9_1}(M,L)$ nor $A_{9_1}(M, L^{-1})$. Hence by Proposition 4.2, there is no 
homomorphism $\pi k_1 \to \pi k_2$ mapping $(m_1, l_1)$ to $(m_2, l_2)$ or 
$(m_2, l_2^{-1})$. \bs

\ni {\bf Acknowledgement.} The second author wishes to thank the
Department of Mathematics of the University of Virginia for its generous
hospitality and the use of its facilities. \bs

\ni {\bf References} \ms

\ni [BZ85] G. Burde and H. Zieschang, Knots, Walter de Gruyter, Berlin,
1985. \ss

\ni [CCGLS94] D. Cooper, M. Culler, H. Gillet, D.D. Long and P.B. Shalen, 
{\it Plane curves associated to character varieties  of knot complements}, {\sl Invent.
Math. \bf 118} (1994), 47--84. \ss

\ni [CL96] D. Cooper and D.D. Long, {\it Remarks on the $A$-polynomial of a knot},
{\sl J. Knot Theory and its Ramif. \bf 5} (1996), 609--628. \ss

\ni [CL98] D. Cooper and D.D. Long, {\it Representation theory and the 
$A$-polynomial of a knot}, {\sl Chaos, Solitons \& Fractals \bf 9} (1998), 749--763. \ss

\ni [CGLS87] M. Culler, C. McA. Gordon, J. Luecke and P.B. Shalen, {\it Dehn surgery on knots}, {\sl Ann. Math. \bf 125} (1987), 237--300.\ss

%\ni [DR99] C. Delman and R. Roberts, {\it Alternating knots satisfy strong property %P}, {\sl Comment. Math. Helv. \bf 74} (1999), 376--397. \ss

\ni [EKT03] S. Eliahou, L. Kauffman and M. Thistlethwaite, {\it Infinite
families of links with trivial Jones polynomial}, {\sl Topology \bf 42}
(2003), 155-169. \ss 

\ni [G83] Gabai, David, {\it Foliations and the topology of $3$-manifolds}, {\sl J.\ Diff.\ Geometry\ \bf 18} (1983), 445--503. \ss

\ni [G-S97] C. Goodman-Strauss, {\it On composite twisted unknots},
{\sl Trans.\ Amer.\ Math.\ Soc.\ \bf349} (1997), 4429--4463.\ss

\ni [GL89] C.McA. Gordon and J. Luecke, {\it Knots are determined by their
complements}, {\sl J. Amer. Math. Soc.} {\bf 2} (1989), 371--415. \ss

%\ni [GPV00] M. Goussarov, M. Polyak and O. Viro, {\it Finite type
%invariants of classical and virtual knots}, Topology {\bf 39}
%(2000), 1045--1068. \ss

%\ni [G82] M. Gromov, Volume and Bounded Cohomology, Inst. Hautes
%\'Etudes Sci. Publ. Math. {\bf 56}, 1982, 5--100. \ss
\ni [H77] R. Hartshorne, Algebraic geometry, Springer, Berlin, 1977. \ss

\ni [H64] J. Hempel, {\it A simply-connected 3-manifold is $S^3$ if it is the sum of a solid torus and the complement of a torus knot}, {\sl Proc. Amer. Math. Soc.} {\bf 15} (1964), 1954--1958. \ss

\ni [H87] J. Hempel, {\it Residual finiteness for $3$-manifolds}, in 
Comb. Group Theory and Topology, {\sl Annals of Math. Studies} {\bf 111}
(1987), 379--396. \ss

%\ni [JS76] W. Jaco and P.B. Shalen, {\it Peripheral structure of
%$3$-manfifolds}, Inventiones math. {\bf 38} (1976), 55--87.\ss

\ni [JS79] W.H. Jaco and P.B. Shalen, {\it Seifert fibered spaces in 
$3$-manifolds}, {\sl Mem.\ Amer.\ Math.\ Soc. \bf21} (1979), no. 220.\ss

\ni [J79]  K.  Johannson, {\it Homotopy equivalence of $3$-manifolds with
boundaries}, {\sl LNM \bf 761}, Springer-Verlag, Berlin, 1979. \ss

\ni [JL89] D. Johnson and C. Livingston, {\it Periperally specified
homomorphisms of knot groups}, {\sl Trans.\ 
Amer.\ Math.\ Soc.\ \bf 311} (1989), 135--146.\ss

\ni [K97] L.H. Kauffman, {\it Talks at MSRI meeting in January
1997}, AMS meeting at University of Maryland, College Park in March
1997, Isaac
Newton Institute Lecture in November 1997, Knots in Hellas Meeting
in Delphi, Greece in July 1998, APCTP-NANKAI Symposium on
Yang-Baxter Systems, Non-Linear Models and Applications at Seoul,
Korea in October 1998.\ss

\ni [K99] L.H. Kauffman, {\it Virtual knot theory}, European J.
Comb. {\bf 20} (1999) 663--690.\ss

\ni [K00] L.H. Kauffman, {\it A survey of virtual knot theory}, in
Knots in Hellas '98, ed. by C. McA. Gordon, V.F.R. Jones, L.H.
Kauffman, S. Lambropoulou, J.H. Przytycki, World Scientific,
Singapore 2000, 143--202.\ss

\ni [Ki00] S.-K. Kim, {\it Virtual knot groups and peripheral
structure}, {\sl J. Knot\ Theory\ Ramif.\ \bf9} (2000), 797--812\ss

\ni [K95] R. Kirby,  {\it Problems in Low-Dimensional Topology}, preprint (1995), \break
http://www.math.berkeley.edu/~kirby/. \ss

\ni [KM04] P.B. Kronheimer and T.S. Mrowka, {\it Witten's conjecture and Property P}, preprint, 2004. \ss

\ni [LS77] R.C.  Lyndon  and P.E. Schupp, Combinatorial Group Theory,
Springer-Verlag, Berlin, 1977. \ss

\ni [MKS76] W. Magnus, A. Karrass and D. Solitar, Combinatorial Group
Theory, Second Edition, Dover, NY, 1976.\ss

\ni [M40] A. I. Mal'cev, {\it On isomorphic matrix representations
of  infinite groups}, {\sl Math. Sb.} {\bf 8} (1940), 405--421. \ss

\ni [M60] K. Murasugi, {\it On alternating knots}, {\sl Osaka\ Math.\ J.\
\bf12} (1960), 277--303.\ss
 
 \ni [N02] F. Nagasato, {\it An approach to the $A$-polynomial of $(2,, 2p+1)$-torus
 knots from Frohman-Gelca-Lofaro theory}, preprint, 2002   http://harold.math.kyushu-u.ac.jp/%~fukky\hfil \break /math.html  \ss

\ni [NZ85] W.D. Neumann and D. Zagier, {\it Volumes of hyperbolic three-manifolds},
{\sl Topology\ \bf 24} (1985), 307--322.\ss

\ni [R76] D. Rolfsen, Knots and Links, Publish or Perish, Inc., Berkeley,
CA 1976. \ss

%\ni [R92] Y.W. Rong, {\it Degree one normal maps between geometric
%$3$-manifolds}, {\sl Trans.\ Amer.\ Math.\ Soc.\ \bf 332} (1992),
%411--436. \ss

\ni [S33] H. Seifert, {\it Topologie dreidimensionaler gefaserter R\"aume}, {\sl Acta Math.\ \bf60} (1933), 147--238. \ss

\ni [S34] H. Seifert, {\it \"Uber das Geschlecht von Knoten}, {\sl
Math.\ Annalen\ \bf 110} (1934), 571--592. \ss

\ni [S95] D.S. Silver, {\it Knot invariants from topological entropy}, 
{\sl Top.\ Appl.\ \bf61} (1995), 159--177. \ss

\ni [SW00] D.S. Silver and S.G. Williams, {\it Virtual knot groups},
in Knots in Hellas '98, Proceedings of the International Conference
on Knot Theory and its Ramifications (C. McA. Gordon, V.F.R. Jones,
L.H. Kauffman, L. Lambropoulou and J.H. Przytycki, eds) World Scientific,
Singapore, 2000, 440--451.\ss

\ni [SW02] D.S. Silver and S.G. Williams, {\it Torsion numbers of
augmented groups}, {\sl L'Enseign. \ Math.\ \bf 48} (2002), 317--343.\ss

%\ni [S81] T. Soma, {\it The Gromov invariant of links}, {\sl
%Invent.\ math.\ \bf 64} (1981), 445--454. \ss

\ni [S02] T. Soma, {\it Epimorphism sequences between hyperbolic
$3$-manifold groups}, {\sl Proc.\ Amer.\ Math.\ Soc.\ \bf130}
(2002), 1221--1223. \ss

\ni [T77/83] W. Thurston, {\it The geometry and topology of $3$-manifolds}
(mimeographed notes), Princeton Univ., 1977/78. \ss

%\ni [T88] V. Turaev, {\it The Yang-Baxter equation and invariants of
%links}, Invent. Math. {\bf 92} (1988) 1225--1237. \ss

\ni [W68] F. Waldhausen, {\it On irreducible $3$-manifolds which are
sufficiently large}, {\sl Annals\ Math.\ \bf 87} (1968), 56--88. \ss

\ni [W02] S. Wang, {\it Non-zero degree maps between 3-manifolds}, Proceedings
of ICM, Volume II, Higher Ed. Press, Beijing, 2002, 457--468. \ss

\ni{\it  Department of Mathematics and Statistics,
University of South Alabama Mobile, AL 36688 USA;
silver@jaguar1.usouthal.edu}\ms
\ni{\it
bjwcw@aol.com}
\end